\documentclass[11pt]{article}
\usepackage{amsfonts}
\usepackage{amsthm,amsmath,amssymb,anysize}
\usepackage[numbers,sort&compress]{natbib}
\newtheorem{lemma}{Lemma}[section]
\newtheorem{theorem}[lemma]{Theorem}
\newtheorem{remark}[lemma]{Remark}
\newtheorem{proposition}[lemma]{Proposition}

\marginsize{3.5cm}{3.5cm}{3.6cm}{3cm}
\numberwithin{equation}{section}

\title{\textsf{Associative forms and second cohomologies of  Lie
superalgebras $HO$ and $KO$}}
 \author{\textsc{Jixia Yuan}$^{1,}$\footnote{Correspondence: jxy@hrbnu.edu.cn (J. Yuan),
 wendeliu@ustc.edu.cn (W. Liu)} \footnote{Supported by  the NSF of Heilongjiang
  Province (A200903) and the NSF of Heilongjiang Educational Committee, China
  (11541109)}
 \ \ \ \textsc{Wende Liu}$^{2,}$\footnote{Supported by the NSF of China (10871057)
and the NSF of Heilongjiang
  Province, China (A200802)} \ \ \ \textsc{Wei Bai}$^{2}$
  \\
  \\
  \ \ \textit{$^{1}$Department of Mathematics},
  \textit{Harbin Institute of Technology}\\
  \textit{Harbin 150006, China}
 \\
   \ \ \textit{$^{2}$School of Mathematical Sciences},
  \textit{Harbin Normal University} \\
  \textit{Harbin 150025, China} }
\date{ }

\begin{document}
\maketitle
\begin{quotation}
\noindent\textbf{Abstract}
We consider two families of finite-dimensional simple Lie superalgebras of Cartan type, denoted
by $HO$ and $KO$, over an algebraically closed
field of  characteristic $p>3.$ Using the weight space
decompositions and the principal gradings we first show that neither
$HO$ nor $KO$ possesses a nondegenerate associative form. Then, by
 means of computing the superderivations from the Lie superalgebras in consideration into their dual modules,
 the second cohomology groups with coefficients in the trivial modules are proved  to be vanishing.\\


  \end{quotation}
\setcounter{section}{-1}
\section{Introduction}
The theory of modular Lie superalgebras has undergone a evolution  during the last ten years, especially in
the
classification of classical simple  modular Lie
 superalgebras (for example, see \cite{BGL,BL}) and the structures and representations of  simple modular Lie superalgebras of Cartan type (for example, see
\cite{fzj,lh,lz,lz1,lzw,wz,xz1,xz,Zh1,Zh2,Zh8,Zh9,Zh10,ZF}).
Recently, one can also find the work on the representations of the classical modular
Lie superalgebras. The present paper is interested
in the second cohomology groups of modular Lie superalgebras of
Cartan type. The classical cohomology  vanishing results of Lie
algebras depend on the characteristic of the underlying field. For
simple Lie superalgebras, even in the case of characteristic zero,
the complete reducibility and the cohomology vanishing theorems do
not hold in general. In view of these general facts above, in this
paper we   prove that the cohomology groups vanish for
 two families of finite-dimensional simple modular Lie superalgebras defined by odd differential forms, known as the odd Hamiltonian
  superalgebras and the odd contact
 superalgebras, respectively. That is,
\begin{itemize}
\item \textsc{Main Theorem} (Theorem \ref{t2}):
The second cohomology group  with coefficients in the trivial module vanishes for both the odd Hamiltonian
  superalgebra  and the odd contact
 superalgebra over an algebraically closed field of characteristic $p>3$.
\end{itemize}

  Let us closing this introduction by briefly recalling certain results on the cohomology
  of modular Lie (super)algebras of Cartan type.
S. Chiu \cite{c}   and R. Farnsteiner \cite {f1,f2} determined the
second cohomology groups of modular Lie algebras of Cartan type. Y. Wang and Y.-Z.
Zhang \cite{wz} determine the second cohomology groups of simple
modular Lie superalgebras of Cartan type $H$ and $K$, which possess
 a nondegenerate associative form  and
W.-J. Xie and Y.-Z. Zhang \cite{xz1,xz} determined the second cohomology groups
of simple modular Lie superalgebras of Cartan type $W,$ $S$ and $K$, which
do not possess
a nondegenerate associative form.

\section{Basics}
 Hereafter  $\mathbb{F}$ is an algebraically closed filed of characteristic
 $p>3,$  $\mathbb{Z}_2:= \{\bar{0},\bar{1}\}$ is the additive group
of order $2$, $\mathbb{Z},$ $\mathbb{N}$ and $\mathbb{N}_{0}$ are
the sets of integers, nonnegative integers and positive integers,
respectively. Throughout this paper $L=L_{\bar{0}}\oplus L_{\bar{1}}=\oplus_{i=-r}^{q}L_{i}$ is a finite
dimensional $\mathbb{Z}$-graded   Lie superalgebra over
$\mathbb{F}$.   Here we should mention
that once the symbol $\mathrm{p}(x)$ appears in an expression, it
 implies that $x$ is  a $\mathbb{Z}_2$-homogeneous element of parity $\mathrm{p}(x)$.
  Write  $L^{\ast}$ for the dual module
of $L$ and $U(L)$ for the enveloping algebra of $L$.
 Then $L^{\ast}$ inherits a $\mathbb{Z}$-graded  $U(L)$-module structure:
\begin{eqnarray*}
&&L^{\ast}=\oplus_{i=-q}^{r}(L^{\ast})_{i};\\
&&(u\cdot f)=(-1)^{\mathrm{p}(u)\mathrm{p}(f)}f\circ\Theta (u)|_{L}
\quad  \mbox{for all}\; u\in U(L),\;f\in L^{\ast},
\end{eqnarray*}
where
$\Theta $
 is the principal anti-automorphism of $U(L)$, that is, the linear even mapping satisfying that
\begin{eqnarray*}
&&\Theta(1)=1;\\
&&\Theta (x)=-x \quad  \mbox{for all}\; x\in L;\\
&&\Theta (uv)=(-1)^{\mathrm{p}(u)\mathrm{p}(v)}\Theta (v)\Theta (u) \quad  \mbox{for all}\;  u,v\in U(L).
\end{eqnarray*}
A superderivation from $L$ into $L$-module $L^{\ast}$ is by definition
a linear mapping $\psi : L\longrightarrow L^{\ast}$ such that
$$\psi([x,y])=(-1)^{\mathrm{p}(\psi)\mathrm{p}(x)}x\cdot \psi(y)
-(-1)^{(\mathrm{p}(\psi)+\mathrm{p}(x))\mathrm{p}(y)}y\cdot \psi(x)\quad \mbox{for all $x,y\in L;$}$$
  $\psi$
 is said to be
inner if there is
some $f\in L^{\ast}  $ such that
$$
\psi(x)=(-1)^{\mathrm{p}(x)\mathrm{p}(f)}x\cdot f\quad \mbox{for all}\; x\in L;
$$
 $\psi$ is said to be skew
if
$$
\psi(x)(y)=-(-1)^{\mathrm{p}(x)\mathrm{p}(y)}\psi(y)(x)\quad \mbox{for all}\;x,y\in L.
$$

Let $H$ be a nilpotent subalgebra of $L_{0}\cap L_{\bar{0}}$ with
weight space decomposition: $L=\oplus_{\alpha\in
\Delta}L_{(\alpha)}.$ In this paper we shall write $\theta$ for zero
weight of $L$. Viewing $L_{i}$ as $H$-module by the adjoint representation,
one can find subsets $\Delta_{i}\subset\Delta$ such that
$L_{i}=\oplus_{\alpha\in \Delta_{i}}L_{i}\cap L_{(\alpha)}.$ Thus
$L$ has a
 $\mathbb{Z}\times \mathrm{Map}(H,\mathbb{F})$-grading structure, which induces a
 $\mathbb{Z}\times \mathrm{Map}(H,\mathbb{F})$-grading structure on  the dual module $L^{\ast}$.

Let $\mathrm{Der}_{\mathbb{F}}(L,L^{\ast})$  be the space of
superderivations from $L$ into $L^{\ast}$ and
$\mathrm{Inn}_{\mathbb{F}}(L,L^{\ast})$  the subspace consisting of
inner superderivations. Then $\mathrm{Der}_{\mathbb{F}}(L,L^{\ast})$
inherits a $\mathbb{Z}\times \mathrm{Map}(H,\mathbb{F})$-grading
from $L$ and $L^{\ast}$ in the usual way.



\section{Associative Forms and Weight Space Decompositions}
 Let $m$ be a  positive integer and suppose that we are given   two $m$-tuples of positive integers,  $ \underline{t}:=\left(
t_{1},t_{2},\ldots,t_m\right) $ and $\pi:=\left( \pi _1,\pi
_2,\ldots,\pi _m\right), $ where $\pi _i:=p^{t_i}-1.$ Let
$\mathcal{O}(m;\underline{t})$ be the divided power algebra over
$\mathbb{F}$ with basis $\{ x^{( \alpha ) }\mid \alpha \in
\mathbb{A}(m;\underline{t})\},$  where
$\mathbb{A}(m;\underline{t}):=\left\{ \alpha \in
\mathbb{N}^m\mid\alpha _i\leq \pi _i \right\}$. For $\varepsilon
_i:=( \delta_{i1}, \ldots,\delta _{im}),$ we abbreviate $x^{(
\varepsilon _i)}$ to $x_i,$ $i=1,\ldots,m.$ Let $\Lambda (n)$ be the
exterior superalgebra over $\mathbb{F}$ with $n$ variables
$x_{m+1},\ldots ,x_{m+n}.$ The tensor product
$\mathcal{O}(m,n;\underline{t})
:=\mathcal{O}(m;\underline{t})\otimes_{\mathbb{F}} \Lambda(n)$
 is an associative super-commutative
superalgebra with a $\mathbb{Z}_2$-grading structure induced by the
trivial $\mathbb{Z}_2$-grading of $\mathcal{O}(m;\underline{t})$ and
the standard $\mathbb{Z}_2$-grading of $\Lambda (n).$   For $g\in
\mathcal{O}(m,\underline{t}),$ $f\in \Lambda(n),$   write
$gf $ for $ g\otimes f$. Note that $ x^{(\alpha) }x^{(\beta)
}=\binom{\alpha +\beta }{ \alpha}
 x^{( \alpha +\beta)}$ for $\alpha,\beta \in
\mathbb{N}^m, $ where $\binom{ \alpha +\beta}
{\alpha}:=\prod_{i=1}^m\binom{ \alpha _i+\beta _i}{ \alpha _i}.$ Let
$$\mathbb{B}(n):=\left\{ \langle i_1,i_2,\ldots
,i_k\rangle \mid m+1\leq i_1<i_2<\cdots <i_k\leq m+n;\; 0\leq k\leq
n\right\}.$$
 For
 $u:=\langle i_1,i_2,\ldots
,i_k\rangle \in \mathbb{B}(n),$ set $ |u| :=k$
 and write $x^u:=x_{i_1}x_{i_2}\cdots x_{i_k}.$
  Notice that we also  denote the set $\{i_1,i_2,\ldots
,i_k\}$ by $u$ itself. For $u\in \mathbb{B}(n+1),$ put

\[\delta_{i\in u}:= \left\{\begin{array}{ll}
0 &\mbox{if}\; i\in u\\
1 &\mbox{if}\; i\notin u,
\end{array}\right.\quad
\| u\|:= \left\{\begin{array}{ll} |u|+1 &\mbox{if}\; 2n+1\in u
\\|u| &\mbox{if}\; 2n+1\notin u.
\end{array}\right.\]
For $u,v\in \mathbb{B}(n) $ with $u\cap v=\emptyset,$
    define $u+v$ to be $w\in \mathbb{B}(n)$ such that
   and $ w = u \cup v.$ If $v\subset u,$ define $u-v$ to be
   $w\in \mathbb{B}(n)$ such that $w=u\setminus v.$
 Note that
$\mathcal{O}(m,n;\underline{t})$ has a standard $\mathbb{F}$-basis
$\{ x^{( \alpha ) }x^u\mid(\alpha,u) \in
\mathbb{A}(m;\underline{t})\times \mathbb{B}(n)\}. $ Let
$\partial_r$ be the superderivation of
$\mathcal{O}(m,n;\underline{t}) $ such that
$\partial_{r}(x^{(\alpha)})=x^{(\alpha-\varepsilon_{r})}$ for $r\in
\overline{1,m}$ and $\partial_{r}(x_{s})=\delta_{rs}$ for $r,s\in
\overline{m+1,m+n}.$ The  generalized Witt superalgebra  $
W\left(m,n;\underline{t}\right)$ is spanned by all $f_r \partial_r,$
where $ f_r\in \mathcal{O}(m,n;\underline{t}),$ $r\in
\overline{1,m+n}.$ Note that $W(m,n;\underline{t} ) $ is a free $
\mathcal{O} \left( m,n;\underline{t}\right)$-module with basis $ \{
\partial_r\mid r\in
 \overline{1,m+n} \}.$ In particular, $W ( m,n;\underline{t} ) $ has a
 standard $\mathbb F$-basis
 $\{x^{(\alpha)}x^{u}\partial_r\mid (\alpha,u,r)\in
\mathbb{A}(m;\underline{t})\times\mathbb{B}(n)\times
\overline{1,m+n}\}.$
  Put $\xi:=|\pi|+n.$ Recall the  standard $\mathbb{Z}$-grading,
$ \mathcal{O}(m,n;\underline{t})=\oplus_{i=
0}^{\xi}\mathcal{O}(m,n;\underline{t})_{\mathbf{s},i},$ where
$$\mathcal{O}(m,n;\underline{t})_{\mathbf{s},i}:={\rm
span}_{\mathbb{F}}\{x^{(\alpha)}x^u\mid |\alpha|+|u|=i, \alpha\in
\mathbb{A}(m;\underline{t}), u\in \mathbb{B}(n)\}.
$$ It induces naturally the standard grading $
W(m,n;\underline{t})=\oplus_{i=-1}^{\xi-1}W(m,n;\underline{t})_{\mathbf{s},i},$ where
$$W(m,n;\underline{t})_{\mathbf{s},i}:={\rm span}_{\mathbb{F}}\{f\partial_j\mid
f\in{\mathcal{O}(m,n;\underline{t})_{\mathbf{s},i+1}}, j\in
\overline{1,m+n}\}.$$ The standard gradings of $
\mathcal{O}(m,n;\underline{t})$ and $W(m,n;\underline{t})$ are also said to be of
type $(1,\ldots,1\mid 1,\ldots,1).$

We  also use the principal grading
$\mathcal{O}(n,n+1;\underline{t})=\oplus_{i=
0}^{\xi+2}\mathcal{O}(n,n+1;\underline{t})_{\mathbf{p},i},$ where
$$
 \mathcal{O}(n,n+1;\underline{t})_{\mathbf{p},i}:={\rm
span}_\mathbb{F}\{x^{(\alpha)}x^u\mid |\alpha|+\|u\|=i,\,\alpha\in
\mathbb{A}(n;\underline{t}),\, u\in \mathbb{B}(n+1)\},$$ and the
principal grading $W(n,n+1;\underline{t})=\oplus_{i=-2}^{\xi+1}
W(n,n+1;\underline{t})_{\mathbf{p},i},$ where
$$ W(n,n+1;\underline{t})_{\mathbf{p},i}:={\rm
span}_{\mathbb{F}}\{f\partial_j\mid
f\in{\mathcal{O}_{\mathbf{p},i+1+\delta_{j,2n+1}}},j\in
\overline{1,2n+1}\}$$(cf. \cite{k2}).  The principal gradings of
$\mathcal{O}(n,n+1;\underline{t})$ and $W(n,n+1;\underline{t})$ are also said to be
of type $(1,\ldots,1\mid 1,\ldots,1,2). $

Define the linear operator
  $\mathrm{T_H}:\mathcal{O}(n,n;\underline{t})\rightarrow
 W(n,n;\underline{t})$ such that
 \[
 \mathrm{T_H}(a):=\sum_{i\in \overline{1,2n}}(-1)^{\mathrm{p}(\partial_i)\mathrm{p}(a)}\partial_i(a)\partial_{i'}
  \quad\mbox{for}\  a\in
 \mathcal{O}(n,n;\underline{t}),
 \]
where \[ i':= \left\{\begin{array}{lll} 2n+1 &\mbox{if}\;
i=0\\
i+n &\mbox{if}\; i\in\overline{1, n}
\\i-n &\mbox{if}\; i\in \overline{n+1, 2n}.
\end{array}\right.\]
 Note that $\mathrm{T_H}$ is  odd and that
\begin{equation*}\label{hee1.1}
   [\mathrm{T_H}(a),\mathrm{T_H}(b)]=\mathrm{T_H}(\mathrm{T_H}(a)(b)) \quad \mbox{for all}\,
   a,b\in\mathcal{O}(n,n;\underline{t}).
\end{equation*}
Then
\[HO(n,n;\underline{t}):=\{\mathrm{T_H}(a)\,|\,a\in\mathcal{O}(n,n;\underline{t})\}
\]   is a finite dimensional simple modular Lie superalgebra, called the
odd Hamiltonian superalgebra (cf. \cite{k2,lzw}).  The  standard  $\mathbb{Z}$-grading is
listed below:
$$
HO(n,n;\underline{t}):=\oplus_{i=-1}^{\xi-2}HO(n,n;\underline{t})_{i},
$$
where
$$
HO(n,n;\underline{t})_{i}:=HO(n,n;\underline{t})\cap
W(n,n;\underline{t})_{\mathbf{s},i}.
$$

Recall the  finite dimensional odd contact superalgebra,  which is a simple   Lie
superalgebra contained in $W(n,n+1;\underline{t}),$ defined as
follows (cf. \cite{fzj,k2}):
$$
KO(n,n+1;\underline{t}):=\{\mathrm{T_{K}}(a)\mid a\in
\mathcal{O}(n,n;\underline{t})\},
$$
where
$$\mathrm{T_{K}}(a):=\mathrm{\mathrm{T_{H}}}(a)+(-1)^{\mathrm{p}(a)}\partial_{2n+1}(a)\mathfrak{D}+ (\mathfrak{D}(a)-2a )\partial_{2n+1}$$
and
\[\mathfrak{D}:=\sum _{i=1}^{2n}x_{i}\partial_{i},\quad
\quad
\mu(i):=\left\{\begin{array}{ll}\bar{0}  &\mbox{if}\; i\in \overline{1,n}\\
\bar{1}&\mbox{if}\; i\in\overline{n+1,2n+1}.
\end{array}\right. \]
We have the following formula (see \cite{fzj,k2}): for $a,b\in \mathcal{O}(n,n+1;\underline{t})$,
\begin{equation}\label{liuee1}
[\mathrm{T_{K}}(a),
\mathrm{T_{K}}(b)]=\mathrm{T_{K}}(\mathrm{T_{K}}(a)(b)-(-1)^{\mathrm{p}(a)}2\partial_{2n+1}(a)b).
\end{equation}
  The principal
$\mathbb{Z}$-graded is   listed below:
$$
KO(n,n+1;\underline{t}):=\oplus_{i=-2}^{\xi}KO(n,n+1;\underline{t})_{i},
$$
where
$$
KO(n,n+1;\underline{t})_{i}:=KO(n,n+1;\underline{t})\cap
W(n,n+1;\underline{t})_{\mathbf{p},i}.
$$
In particular,
$$
KO(n,n+1;\underline{t})_{-2}=\mathbb{F}\mathrm{T_{K}}(1),\quad
KO(n,n+1;\underline{t})_{\xi}=\mathbb{F}\mathrm{T_{K}}(x^{(\pi)}x^{\langle
n+1,\ldots,2n+1\rangle}).
$$ For simplicity,  we  write $HO$
 and $KO$ for $HO(n,n;\underline{t})$ and
 $KO(n,n+1;\underline{t}),$ respectively.

For  $i\in
\overline{1,n} $, set
$$H_{HO}:=\sum_{i\in
\overline{1,n}}\mathbb{F}\mathrm{T_H}(x_{i}x_{i'}),\quad H_{KO}:=\sum_{i\in
\overline{1,n}}\mathbb{F}\mathrm{T_{K}}(x_{i}x_{i'}).$$
Let $X:=HO$ or $ KO.$ Then $H_{X}$ is a  torus of $X_{\bar{0}}\bigcap X_0$, known as standard.
 We consider the weight space decomposition of $X$ with respect to the torus $H_{X},$
 $$X=\oplus_{\gamma\in \Delta_{X}}X_{\gamma}.$$
For $\alpha\in \mathbb{A}(n;\underline{t}) $ and $u\in \mathbb{B}(n),$ we define $\alpha+u$ to
be the linear function on $H_{X} $,
\begin{eqnarray*}
 \alpha+u: H_{X}\longrightarrow \mathbb{F}
\end{eqnarray*}
such that for  $i\in
\overline{1,n}, $
\begin{eqnarray*}
&&(\alpha+u)(\mathrm{T_H}(x_{i}x_{i'})):=\delta_{i'\in
u}-\alpha_{i}\quad \mbox{if} \; X=HO,\\
&&
 (\alpha+u)(\mathrm{\mathrm{T_{K}}}(x_{i}x_{i'})):=\delta_{i'\in
u}-\alpha_{i}\quad \mbox{if} \; X=KO.
\end{eqnarray*}
Since
\begin{eqnarray*}[\mathrm{T_H}(x_{i}x_{i'}),\mathrm{T_H}(x^{(\alpha)}x^{u})]=
(\delta_{i'\in u}-\alpha_{i})\mathrm{T_H}(x^{(\alpha)}x^{u}),
\end{eqnarray*}
one sees that $\mathrm{T_{H}}(x^{(\alpha)}x^{u})\in HO$ is a weight
vector belonging to the weight $\alpha+u.$ Similarly,
$\mathrm{T_{K}}(x^{(\alpha)}x^{u})$ and
$\mathrm{T_{K}}(x^{(\alpha)}x^{u}x_{2n+1})\in KO$ is a weight
vector belonging to the weight $\alpha+u.$

For $\alpha,\beta\in\mathbb{A}(n;\underline{t}),$ write
$\alpha\equiv\beta\pmod{p}$ provided that
$\alpha_{i}\equiv\beta_{i}\pmod{p} $ for all $i\in \overline{1,n}.$ We have:
\begin{proposition}\label{l8}
For $\alpha\in \mathbb{A}(n;\underline{t}) $ and $u\in
\mathbb{B}(n),$ we have
\begin{eqnarray*}
HO_{(\alpha+u)}=\sum_{\beta\equiv \alpha \pmod{p}, v_{1}\subset u,\;
v_{2}\subset \overline{n+1,2n}\backslash u,\atop\beta- \sum_{j'\in
v_{1} }\varepsilon_{j}+\sum_{j'\in v_{2}}\varepsilon_{j}\in
\mathbb{A}(n;\underline{t}) }\mathbb{F}\mathrm{T_{H}}(x^{(\beta-
\sum_{j'\in v_{1} }\varepsilon_{j}+\sum_{j'\in
v_{2}}\varepsilon_{j})}x^{u-v_{1}+v_{2}}),
\end{eqnarray*}
\begin{eqnarray*}
KO_{(\alpha+u)}=\sum_{\beta\equiv \alpha \pmod{p}, v_{1}\subset u,\;
v_{2}\subset \overline{n+1,2n+1}\backslash u\atop \beta-\sum_{j'\in
v_{1} }\varepsilon_{j}+\sum_{2n+1\neq j'\in v_{2}}\varepsilon_{j}\in
\mathbb{A}(n;\underline{t})}\mathbb{F}\mathrm{\mathrm{T_{K}}}(x^{(\beta-\sum_{j'\in
v_{1} }\varepsilon_{j}+\sum_{2n+1\neq j'\in
v_{2}}\varepsilon_{j})}x^{u-v_{1}+v_{2}}).
\end{eqnarray*}
\end{proposition}

\begin{proof}
(i) For $HO$, the inclusion $``\supset"$ is straightforward. To show
the converse, let
 $\mathrm{T_H}(x^{(\beta)}x^{v})\in HO_{(\alpha+u)},$ where
  $\beta\in \mathbb{A}(n,\underline{t}),$ $v\in \mathbb{B}(n).$ Then, for  $i\in \overline{1,n},$ one has
\begin{eqnarray*}
(\delta_{i'\in v}-\beta_{i})\mathrm{T_H}(x^{(\beta)}x^{v})
&=&(\beta+v)(\mathrm{T_H}(x_{i}x_{i'}))\mathrm{T_H}(x^{(\beta)}x^{v})\\
&=&[\mathrm{T_H}(x_{i}x_{i'}),\mathrm{T_H}(x^{(\beta)}x^{v})]
\\
&=&(\alpha+u)(\mathrm{T_H}(x_{i}x_{i'}))\mathrm{T_H}(x^{(\beta)}x^{v})
\\&=&(\delta_{i'\in
u}-\alpha_{i})\mathrm{T_H}(x^{(\beta)}x^{v}).
\end{eqnarray*}
 Consequently,  $\delta_{i'\in
v}-\beta_{i}=\delta_{i'\in u}-\alpha_{i}$ holds in $\mathbb{F}$, that is,
\[  \alpha_{i}-\beta_{i}=\delta_{i'\in
u}-\delta_{i'\in v}= \left\{\begin{array}{lll} 1\pmod{p}
&\mbox{if}\; i'\in u \;\mbox{and}\; i'\notin v
\\-1\pmod{p} &\mbox{if}\; i'\notin u \;\mbox{and}\; i'\in v
\\0\pmod{p} &\mbox{if}\; i'\in u \;\mbox{and}\; i'\in v
\\0\pmod{p} &\mbox{if}\; i'\notin u \;\mbox{and}\; i'\notin v.
\end{array}\right.\]
The assertion follows.

(ii) In (i), replacing $HO$ and $\mathrm{T_H}$ by $KO$ and $\mathrm{T_{K}}$, respectively, one checks that all the arguments above still
 hold.
Thus, for $KO$ it can be proved in the same manner.
\end{proof}

Write $\omega_{0}:=\langle n+1,\ldots,2n\rangle$ and
$\omega_{1}:=\langle n+1,\ldots,2n+1 \rangle.$ By Proposition \ref{l8}, we list the following for later use:
\begin{eqnarray}
&&KO_{\theta}=\sum_{\alpha\equiv0\pmod{p},\;u\in
\mathbb{B}(n+1)}\mathbb{F}\mathrm{T_{K}}(x^{(\alpha+\sum_{2n+1\neq j'\in u}\varepsilon_{j})}x^{u}),\label{e4}\\
&&\Delta_{HO,0}=\Delta_{KO,0}=\big\{\theta, 2\varepsilon_{i},
\varepsilon_{i}+\varepsilon_{j}, \varepsilon_{i}+\langle j'\rangle,
\langle i',j'\rangle\mid i\neq j\in \overline{1,n}\big\},\label{e5}\\
&&\Delta_{HO,1}=\Delta_{KO,1}=\big\{\varepsilon_{i}, \langle i'\rangle,
2\varepsilon_{i}+\varepsilon_{j}, 2\varepsilon_{i}+\langle j\rangle,
\varepsilon_{i}+\varepsilon_{j}+\varepsilon_{k},
\varepsilon_{i}+\varepsilon_{j}+\langle k'\rangle,\nonumber\\
&&~~~~~~~~~~~~~~~~~~~~~~~~~~~ 3\varepsilon_{i},
\varepsilon_{i}+\langle j',k'\rangle, \langle i',j',k'\rangle\mid
\,\mbox{distinct}\, i, j, k\in
\overline{1,n}\big\},\label{e6}\\
&&\Delta_{HO,2}=\Delta_{KO,2}=\Big\{\theta,\alpha+u,\beta+v\mid
|\alpha|+|u|=4,
\sum_{i=1}^{n}\partial_{i}\partial_{i'}(x^{(\alpha)}x^{u})=0,\nonumber\\
&&~~~~~~~~~~~~~~~~~~~~~~~~~~~~~~~~~~~~~~~~~~~~~~~~|\beta|+|v|=2,
\sum_{i=1}^{n}\partial_{i}\partial_{i'}(x^{(\beta)}x^{v})=0\Big\},\label{e7}\\
&&\Delta_{HO,\xi-2}=\Delta_{KO,\xi}=\{\pi+\omega_{0}\}.\label{e8}
\end{eqnarray}

 As in the Lie algebra case, one may easily prove the following
\begin{lemma}\label{l6}
Let $L=\oplus_{i=-r}^{q}L_{i} $ be a finite dimensional simple
$\mathbb{Z}$-graded Lie superalgebra and $H$ a
 nilpotent subalgebra of
$L_{\bar{0}}\cap L_{0}$ with weight space decomposition
$L=\oplus_{\alpha\in \Delta}L_{(\alpha)}.$ If $L$ has a
nondegenerate associative form, then
  $$\dim_{\mathbb{F}}L_{k}=\dim_{\mathbb{F}}L_{q-r-k} \quad\mbox{and}\quad\dim_{\mathbb{F}}L_{k}\cap
L_{\gamma}=\dim_{\mathbb{F}}L_{q-r-k}\cap L_{-\gamma} $$ for all
$-r\leq k \leq q $ and $\gamma\in \Delta$.
\end{lemma}

\begin{theorem}\label{t1} Let $X:=HO$ or $KO$. Then
$X$ does not possess any nondegenerate associative form.
\end{theorem}
\begin{proof}
For $X=HO,$ we have $X=\sum_{i=-1}^{\xi-2}X_{i}.$ Note that
\begin{eqnarray*}
 X_{\xi-2}=\mathrm{span}_{\mathbb{F}}{\mathrm{T_{H}}(x^{(\pi)}x^{\omega_{0}})},\quad
 X_{-1}=\mathrm{span}_{\mathbb{F}}\{\mathrm{T_{H}}(x_{i})\mid i\in
\overline{1,2n}\}.
 \end{eqnarray*} This shows that
 $$1=\dim
X_{\xi-2}\neq \dim X_{-1}=2n.$$

For $X=KO,$ we have $X=\sum_{i=-2}^{\xi}X_{i}$. Consider the weight
space decomposition with respect to the torus $H_{{X}} $
defined above,
$$
 X=\oplus_{\alpha\in \Delta_{X}}X_{(\alpha)}.
$$
From (\ref{e4}) we have
$$1=\dim X_{-2}\cap X_{\theta}\neq \dim X_{\xi}\cap
X_{\theta}=0.$$
 Now our conclusion follows from
Lemma \ref{l6}.
\end{proof}
\section{Second Cohomology Groups}

As in Lie algebra case \cite{f2}, if a finite dimensional simple Lie algebra $L$
 does not possess any nondegenerate associative form then $H^{2}(L,\mathbb{F})\cong H^{1}(L,L^{\ast})$ and
all the  $\mathbb{Z}_{2}$-homogeneous superderivations from $L$ into $L^*$ are skew. We thereby obtain by
 Lemma \ref{l6} and Theorem \ref{t1} the following
\begin{theorem}\label{liu1631}
 $H^{2}(X,\mathbb{F})$ is isomorphic to $H^{1}(X,X^{\ast})$ and  $\mathrm{Der}(X,X^{\ast})$ consists of
 skew superderivations for $X:=HO$ or $KO$.
\end{theorem}

In view of this theorem, in order to determine the second cohomology
groups, it is enough to compute $H^{1}(X,X^{\ast})$. For this
purpose, we want to establish a reduction lemma relative to the
superderivations from $X$ into $X^{*} $ (see Lemma \ref{l7}). Before
doing that, we first recall a general result (Lemma \ref{l2}),
 which is
analogous to the Lie algebra case \cite{f2}.

Let $L$ be a finite
dimensional $\mathbb{Z}$-graded  Lie superalgebra over $\mathbb{F}$, $L=\oplus_{i=-r}^{q}L_{i}.$
Write $L^{-}:=\oplus_{i=-r}^{-1} L_{i}$ and $L^{+}:=\oplus_{i=1}^{q}L_{i}.$
Recall   the
canonical mapping
$$\Phi_{1}: H^{1}(L,L^{\ast})\longrightarrow
H^{1}(L^{-},L^{\ast}) $$
which is naturally  induced by the restriction mapping
$$
\mathrm{Der}_{\mathbb{F}}(L,L^{\ast})\longrightarrow
\mathrm{Der}_{\mathbb{F}}(L^{-},L^{\ast}).
$$
 Write $U(L)^{+}$ for the two-sided ideal of  $U(L)$ generated by $L $ and put
 $$\mathrm{Ann}_{{U(L^{-})}^{+}}(L):=\{u\in U(L^{-})^{+}\mid u\cdot
L=0\}.$$

\begin{lemma} \cite{xz}\label{l2}
Let $V$ be a $\mathbb{Z}_{2}$-graded subspace of $L$. Suppose there
are a $\mathbb{Z}_{2}$-homogeneous basis $\{e_{1},e_{2},\ldots,
e_{n}\}$ of $L^{-}$ and   a $\mathbb{Z}_{2}$-homogeneous basis
$\{v_{1},v_{2},\ldots, v_{m}\}$ of $V$ such that  $\{e^{a}\cdot
v_{j}\mid  1\leq j\leq m, a\in T\}$ is a basis of $L$ over
$\mathbb{F} $ for some  subset $T\subset \mathbb{N}^{n}  $ and that
$\mathrm{Ann}_{{U(L^{-})}^{+}}(L)=
\mathrm{span}_{\mathbb{F}}\{e^{b}\mid b \in\mathbb{N}^{n}\setminus
T\}.$
  Then the following statements hold:

\begin{itemize}
\item[$\mathrm{(1)}$]  Every superderivation $\psi: L\longrightarrow
L^{\ast}$ satisfying the condition that $\mathrm{ker}(\mathrm{ad}e_{i})\subset
\mathrm{ker}\psi(e_{i}) $ for all $i\in \overline{1,n}$ defines an
element of $\mathrm{ker}\Phi_{1},$ that is, $\psi\mid_{L^{-}}\in
\mathrm{Inn}_{\mathbb{F}}(L^{-},L^{\ast}).$

\item[$\mathrm{(2)}$]   Suppose there is $\mu\in \mathbb{N}^{n}$ such that
$T=\{b\in \mathbb{N}^{n}\mid b\leq \mu\}.$ Then
$\mathrm{ker}(\mathrm{ad}e_{i})\subset \mathrm{ker}\psi(e_{i})$ if
and only if $e_{i}^{\mu_{i}}\cdot \psi(e_{i})=0.$
\end{itemize}
\end{lemma}

We  mention a standard result  that for a Lie superalgebra $L$  every homogeneous superderivation
$\psi: L\longrightarrow L^{\ast}$ extends to  one and only one
$U(L)$-module homomorphism $\Psi: U(L)^{+}\longrightarrow L^{\ast}$
so that $\Psi(x)=\psi(x)$ for all $x\in L.$
\begin{lemma}\label{l9}
If $\psi: HO\longrightarrow HO^{\ast}$ is a  superderivation then
$$
\mathrm{T_{H}}(x_{i})\cdot\psi(\mathrm{T_{H}}(x_{i}))=\mathrm{T_{H}}(x_{j})^{\pi_{j}}\cdot\psi(\mathrm{T_{H}}(x_{j}))=0\quad\mbox{for
all }\; i\in \overline{1,n},\; j\in \overline{n+1,2n}.
$$
\end{lemma}
\begin{proof}
Letting $\Psi$ be as indicated above,    we have
\begin{eqnarray*}
&&\mathrm{T_{H}}(x_{i})\cdot\psi(\mathrm{T_{H}}(x_{i}))=
\mathrm{T_{H}}(x_{i})\cdot\Psi(\mathrm{T_{H}}(x_{i}))=\Psi(\mathrm{T_{H}}(x_{i})^{2})=0
\;\;\;\mbox{for all} \;i\in \overline{1,n},\\
&&\mathrm{T_{H}}(x_{j})^{\pi_{j}}\cdot\psi(\mathrm{T_{H}}(x_{j}))=
\mathrm{T_{H}}(x_{j})^{\pi_{j}}\cdot\Psi(\mathrm{T_{H}}(x_{j}))
=\Psi(\mathrm{T_{H}}(x_{j})^{p^{t_{j}}}) \;\mbox{for all}\;\;\; j\in
\overline{n+1,2n}.
\end{eqnarray*}
Let $G\in HO.$ As
$\mathrm{T_{H}}(x_{j})^{p^{t_{j}}}$ lies in the center
 $C(U(HO)^{+}),$ where $j\in \overline{n+1,2n},$ we have
$$
G\cdot
\Psi(\mathrm{T_{H}}(x_{j})^{p^{t_{j}}})=\Psi(G\mathrm{T_{H}}(x_{j})^{p^{t_{j}}})
=\Psi(\mathrm{T_{H}}(x_{j})^{p^{t_{j}}}G)
=\mathrm{T_{H}}(x_{j})^{p^{t_{j}}}\cdot\Psi(G).
$$
Note that
$$
(\mathrm{T_{H}}(x_{j})^{p^{t_{j}}}\cdot\Psi(G))(y)=\pm
\Psi(G)(\mathrm{T_{H}}(x_{j})^{p^{t_{j}}}\cdot y)=0
$$
for every $y\in HO.$ Consequently, $\mathrm{T_{H}}(x_{j})^{\pi_{j}}\cdot\psi(\mathrm{T_{H}}(x_{j}))$ lies in
$$
\{f\in HO^{\ast}\mid HO\cdot f=0\}=\{f\in HO^{\ast}\mid
f([HO,HO])=0\},$$
which is zero, since $HO$ is simple. The proof is complete.
\end{proof}
\begin{lemma}\label{l12}
For $\mu:=(1,\ldots,1,\pi_{1},\ldots,\pi_{n})\in \mathbb{N}^{2n},$
we have
$$\mathrm{ker}(\mathrm{ad}\mathrm{T_{H}}(x_{i}))=\mathrm{T_{H}}(x_{i})^{\mu_{i}}\cdot
HO+\mathbb{F}\mathrm{T_{H}}(x_{i'}),\quad i\in \overline{1,2n}.$$

\end{lemma}
\begin{proof}
The inclusion $``\supset"$ is clear. To show the converse, for
$b:=(b_{1},b_{2},\ldots,b_{2n})\in \mathbb{N}^{2n},$ where $b_{i}=0$
or $1$ for all $i\in \overline{1,n},$ put
\begin{eqnarray*}
&&\mathrm{T_{H}}^{b}:=\mathrm{T_{H}}(x_{1})^{b_{1}}\mathrm{T_{H}}(x_{2})^{b_{2}}\cdots
\mathrm{T_{H}}(x_{2n})^{b_{2n}},\\
&&U(HO_{-1})_{(k)}:=\mathrm{span}_{\mathbb{F}}\Big\{\mathrm{T_{H}}^{b}\mid
 \sum_{i=1}^{2n}b_{i}\leq k\Big\},\\
&&HO_{(k)}:=U(HO_{-1})_{(k)}\cdot\mathrm{T_{H}}(x^{(\pi)}x^{\omega_{0}}).
\end{eqnarray*}
 Then
$HO=\sum_{k=0}^{|\mu|-1}HO_{(k)}.$ By induction one gets
$$\mathrm{ker}(\mathrm{ad}\mathrm{T_{H}}(x_{i}))\cap HO_{(k)}\subset
\mathrm{T_{H}}(x_{i})^{\mu_{i}}\cdot HO,$$ where $0\leq k\leq
|\mu|-2.$ Now we want to show that
$$\mathrm{ker}(\mathrm{ad}\mathrm{T_{H}}(x_{i}))\cap
HO_{(|\mu|-1)}\subset \mathrm{T_{H}}(x_{i})^{\mu_{i}}\cdot
HO+\mathbb{F}\mathrm{T_{H}}(x_{i'}).$$ Take
$$x=\sum_{0\leq
a<\mu}\lambda_{a}\mathrm{T_{H}}^{a}\cdot
\mathrm{T_{H}}(x^{(\pi)}x^{\omega_{0}})\in
\mathrm{ker}(\mathrm{ad}\mathrm{T_{H}}(x_{i}))\cap HO_{(|\mu|-1)}.$$
Then
\begin{eqnarray*}
x&\equiv&\sum_{0\leq
a<\mu,a_{i}=\mu_{i}\atop|a|=|\mu|-1}\lambda_{a}\mathrm{T_{H}}^{a}\cdot
\mathrm{T_{H}}(x^{(\pi)}x^{\omega_{0}})+\lambda_{\mu-\varepsilon_{i}}\mathrm{T_{H}}^{\mu-\varepsilon_{i}}\cdot\mathrm{T_{H}}(x^{(\pi)}x^{\omega_{0}})\\
&\equiv& \sum_{0\leq
a<\mu,a_{i}=\mu_{i}\atop|a|=|\mu|-1}\pm\lambda_{a}
\mathrm{T_{H}}(x_{i})^{\mu_{i}}\mathrm{T_{H}}^{a-\mu_{i}\varepsilon_{i}}\cdot
\mathrm{T_{H}}(x^{(\pi)}x^{\omega_{0}})\pm\lambda_{\mu-\varepsilon_{i}}\mathrm{T_{H}}(x_{i'})\pmod{HO_{(|\mu|-2)}}.
\end{eqnarray*}
Therefore, there is $y\in HO$ such that
$$x-\mathrm{T_{H}}(x_{i})^{\mu_{i}}\cdot
y\pm\lambda_{\mu-\varepsilon_{i}}\mathrm{T_{H}}(x_{i'})\in
HO_{(|\mu|-2)}
$$
and
$$\mathrm{ad}\mathrm{T_{H}}(x_{i})(x-\mathrm{T_{H}}(x_{i})^{\mu_{i}}\cdot
y\pm\lambda_{\mu-\varepsilon_{i}}\mathrm{T_{H}}(x_{i'}))=0.
$$
Then
$$x-\mathrm{T_{H}}(x_{i})^{\mu_{i}}\cdot
y\pm\lambda_{\mu-\varepsilon_{i}}\mathrm{T_{H}}(x_{i'})\in\mathrm{ker}(\mathrm{ad}\mathrm{T_{H}}(x_{i}))\cap
HO_{(|\mu|-2)}\subset \mathrm{T_{H}}(x_{i})^{\mu_{i}}\cdot HO,$$
that is, $x\in \mathrm{T_{H}}(x_{i})^{\mu_{i}}\cdot
HO+\mathbb{F}\mathrm{T_{H}}(x_{i'}).$ The proof is complete.
\end{proof}
Recall that for $X=HO$ or $ KO,$   $H_{X}$ is a torus of $X_{\bar{0}}\bigcap X_0$.
Consider the weight space decomposition of $X$ with respect to $H_{X},$
 $$X=\oplus_{\gamma\in \Delta_{X}}X_{\gamma}.$$
 We would like to mention  a standard fact that every
$\mathrm{Map}(H_{X},\mathbb{F})$-homogeneous nonzero-degree
superderivation from $X$ into $X^{\ast}$ must be inner. So we are concerned with those superderivations of
$\mathrm{Map}(H_{X},\mathbb{F})$-degree $\theta.$

\begin{lemma}\label{l7}
Let $X$ stand for $HO$ or $KO.$ If $\psi: X\longrightarrow X^{\ast}$
is a superderivation of $\mathrm{Map}(H_{X},\mathbb{F})$-degree $\theta,$ then there exists some
$f\in X^{\ast}$ such that
$\psi(x)=(-1)^{\mathrm{d}(x)\mathrm{d}(f)}x\cdot f$ for all $x\in
X^{-}.$
\end{lemma}
\begin{proof} (i) First consider the case $X=HO.$ The general assumption of Lemma
\ref{l2} is valid for
$V:=\mathbb{F}\cdot\mathrm{T_{H}}(x^{(\pi)}x^{\omega_{0}}).$  Put
$\mathrm{T_{H}}^{a}:=\mathrm{T_{H}}(x_{1})^{a_{1}}\mathrm{T_{H}}(x_{2})^{a_{2}}\cdots
\mathrm{T_{H}}(x_{2n})^{a_{2n}}$ for
$a:=(a_{1},a_{2},\ldots,a_{2n})\in \mathbb{N}^{2n},$ where $a_{i}=0$
or $1$ for all $i\in \overline{1,n}.$ Suppose
$a_{i_{1}}=a_{i_{2}}=\cdots = a_{i_{k}}=0,$ $1\leq i_{1}<
i_{2}<\cdots i_{k}\leq n.$ Write
$$
x^{u}:=x_{i_{1}'}x_{i_{2}^{'}}\cdots x_{i_{k}^{'}}, \quad
b:=a_{1'}\varepsilon_{1}+a_{2'}\varepsilon_{2}+\cdots
a_{n'}\varepsilon_{n}.
$$
Since
$$
\mathrm{T_{H}}^{a}\cdot \mathrm{T_{H}}(x^{(\pi)}x^{\omega_{0}})=\pm
\mathrm{T_{H}}(x^{(\pi-b)}x^{u}),
$$
$\{\mathrm{T_{H}}^{a}\cdot
\mathrm{T_{H}}(x^{(\pi)}x^{\omega_{0}})\mid a\in T\} $ is an
$\mathbb{F}$-basis of $HO$ and
$$
\mathrm{Ann}_{U(HO^{-})^{+}}(HO)=\{\mathrm{T_{H}}^{a}\mid a\notin
T\},
$$
where
$$
T:=\{a\in \mathbb{N}^{2n}\mid a< \mu\} ,\quad
\mu:=(1,\ldots,1,\pi_{1},\ldots,\pi_{n})\in\mathbb{N}^{2n}.
$$
 By Lemmas \ref{l2}(1) and \ref{l12}, it suffices to show that
$ \mathrm{T_{H}}(x_{i})^{\mu_{i}}\cdot
HO+\mathbb{F}\mathrm{T_{H}}(x_{i'})\subset
\mathrm{ker}\psi(\mathrm{T_{H}}(x_{i})),$ for all  $i\in
\overline{1,2n}$.
 By Lemma \ref{l9},
\begin{eqnarray*}
0=\mathrm{T_{H}}(x_{i})^{\mu_{i}}\psi(\mathrm{T_{H}}(x_{i}))(HO)
&=&\pm\psi(\mathrm{T_{H}}(x_{i}))(\Theta(\mathrm{T_{H}}(x_{i})^{\mu_{i}})\cdot
HO)\\
&=&\pm\psi(\mathrm{T_{H}}(x_{i}))(\mathrm{T_{H}}(x_{i})^{\mu_{i}}\cdot
HO),
\end{eqnarray*}
that is,  $\mathrm{T_{H}}(x_{i})^{\mu_{i}}\cdot
HO\subset \mathrm{ker}\psi(\mathrm{T_{H}}(x_{i})).$ Since $\psi$ is
a superderivation of degree $\theta,$
 it is clear that
$\mathbb{F}\mathrm{T_{H}}(x_{i'})\subset
\mathrm{ker}\psi(\mathrm{T_{H}}(x_{i})).$

 (ii) Consider the case $X=KO.$
The general assumption of Lemma \ref{l2} is valid for
$$V:=\mathbb{F}\cdot\mathrm{T_{H}}(x^{(\pi)}x^{\omega_{1}}).$$  Put
$$\mathrm{T_{K}}^{a}:=\mathrm{T_{K}}(1)^{a_{0}}\mathrm{T_{K}}(x_{1})^{a_{1}}\mathrm{T_{K}}(x_{2})^{a_{2}}\cdots
\mathrm{T_{K}}(x_{2n})^{a_{2n}}$$ for
$a:=(a_{0},a_{1},a_{2},\ldots,a_{2n})\in \mathbb{N}^{2n+1},$ where $a_{i}=0$
or $1$ for all $i\in \overline{0,n}.$ Suppose
$a_{i_{1}}=a_{i_{2}}=\cdots = a_{i_{k}}=0,$ $0\leq i_{1}<
i_{2}<\cdots i_{k}\leq n.$ Write
$$
x^{u}:=x_{{i_{1}}'}x_{{i_{2}}'}\cdots x_{{i_{k}}'}, \quad
b:=a_{1'}\varepsilon_{1}+a_{2'}\varepsilon_{2}+\cdots
a_{n'}\varepsilon_{n}.
$$
Since
$$
\mathrm{T_{K}}^{a}\cdot \mathrm{T_{K}}(x^{(\pi)}x^{\omega_{1}})=\pm\lambda
\mathrm{T_{K}}(x^{(\pi-b)}x^{u}), \;\mbox{where}\; \lambda:=1 \;\mbox{or}\; 2,
$$
one sees that $ \{\mathrm{T_{K}}^{a}\cdot
\mathrm{T_{K}}(x^{(\pi)}x^{\omega_{1}})\mid a\in T\}$ is
an $\mathbb{F}$-basis of $KO$ and
$$
\mathrm{Ann}_{U(KO^{-})^{+}}(KO)=\{\mathrm{T_{K}}^{a}\mid a\notin
T\},
$$
where
$$
T:=\{a\in \mathbb{N}^{2n+1}\mid a\leq \mu\} ,\quad
\mu:=(1,\ldots,1,\pi_{1},\ldots,\pi_{n})\in\mathbb{N}^{2n+1}.
$$
As in Lemma
\ref{l9}, one easily  get,  for $i\in \overline{1,2n},$
$$
\mathrm{T_{K}}(1)\psi(\mathrm{T_{K}}(1))=\mathrm{T_{K}}(x_{i})^{\mu_{i+1}}\psi(\mathrm{T_{K}}(x_{i}))=0.
$$  Then the desired result  follows from Lemma \ref{l2}(1)
and (2).
\end{proof}
As before, suppose $L=\oplus_{i=-r}^{q}L_{i} $ is a finite
dimensional simple $\mathbb{Z}$-graded Lie superalgebra and $H$ a
 nilpotent subalgebra of
$L_{\bar{0}}\cap L_{0}$ with weight space decomposition
$L=\oplus_{\alpha\in \Delta}L_{(\alpha)}.$  We consider the
subalgebra
$$
 M(L):=[L^{+},L^{+}].
$$
Note that $M(L)$ is a graded subalgebra of $L$ on which $H$
operates. Hence for $h\geq 1$ there is $\phi_{h}\subset \Delta_{h}$
such that $L_{h}=M(L)_{h}+\oplus_{\alpha\in
\phi_{h}}L_{(\alpha)}\cap L_{h}.$

\begin{lemma}\label{l10}
Suppose $h\geq 3$. Then
\begin{equation}\label{eqlt1636}
HO_{h}=M(HO)_{h}+\sum_{\alpha_{i}\equiv 0,1(\mathrm{mod}{p})\; \forall i\in \overline{1,n}
\atop u\in
\mathbb{B}(n),|\alpha|+|u|-2=h}\mathbb{F}\mathrm{T_H}(x^{(\alpha)}x^{u}).
\end{equation}
\end{lemma}
\begin{proof}
The inclusion $``\supset"$ is clear. Let us consider the converse
inclusion. Let $\mathrm{T_{H}}(x^{(\alpha)}x^{u})\in HO_{h},$ where
$h\geq3,$ $\alpha\in \mathbb{A}(n,\underline{t})$ and $u\in
\mathbb{B}(n).$ If $\alpha_{i}\equiv 0,1(\mathrm{mod}{p})$ for all $i\in \overline{1,n}$,
then $\mathrm{T_{H}}(x^{(\alpha)}x^{u})$ lies in the second summand of (\ref{eqlt1636}). Thereby we
suppose  $\alpha_{j}\not\equiv 0,1\pmod{p}$ for some $j\in \overline{1,n}$. We distinguish two cases:\\

\noindent\textit{Case 1.} Suppose $j'\notin u$. We have
\begin{eqnarray*}
0\neq{\alpha\choose
2\varepsilon_{j}}\mathrm{T_H}(x^{(\alpha)}x^{u})=-
[\mathrm{T_H}(x^{(2\varepsilon_{j})}x_{j'}),
\mathrm{T_H}(x^{(\alpha-\varepsilon_{j})}x^{u})]\in M(HO)_{h}.
\end{eqnarray*}

\noindent\textit{Case 2.} Suppose $j'\in u.$ We have
\begin{eqnarray}\label{e1}
\frac{\alpha_{j}(3-\alpha_{j})}{2}\mathrm{T_H}(x^{(\alpha)}x^{u})=\left[\mathrm{T_H}(x^{(2\varepsilon_{j})}x_{j'}),
\mathrm{T_H}(x^{(\alpha-\varepsilon_{j})}x^{u})\right]\in
M(HO)_{h}
\end{eqnarray}
and
\begin{eqnarray}\label{e2}
\frac{\alpha_{j}(\alpha_{j}-1)(5-\alpha_{j})}{6}\mathrm{T_H}(x^{(\alpha)}x^{u})
=\left[\mathrm{T_H}(x^{(3\varepsilon_{j})}x_{j'}),
\mathrm{T_H}(x^{(\alpha-2\varepsilon_{j})}x^{u})\right] \in
M(HO)_{h}.
\end{eqnarray}
Since $\alpha_{j}\not\equiv0,1\pmod{p},$
$\frac{\alpha_{j}(3-\alpha_{j})}{2}$ and
$\frac{\alpha_{j}(\alpha_{j}-1)(5-\alpha_{j})}{6}$ cannot be all
zero modulo $p$. Thus, it follows from  (\ref{e1}) or (\ref{e2})
that $\mathrm{T_H}(x^{(\alpha)}x^{u})\in M(HO)_{h}.$

Summarizing, the proof is complete.
\end{proof}
\begin{remark}\label{r1}
Let $\mathcal{K}(n,\underline{t})$ be the subspace spanned by the
elements $\mathrm{T_{K}}(a) $ with $a\in \oplus_{i\geq 1}\mathcal
{O}(n,n;\underline{t})_{\mathbf{s},i}.$ If $a, b\in \oplus_{i\geq
1}\mathcal {O}(n,n;\underline{t})_{\mathbf{s},i},$   by
(\ref{liuee1}), we have
\begin{equation*}\label{hee1.1}
   [\mathrm{T_{K}}(a),\mathrm{T_{K}}(b)]=\mathrm{T_{K}}(\mathrm{T_{H}}(a)(b)).
\end{equation*}
Note that
\begin{equation*}\label{hee1.1}
   [\mathrm{T_{H}}(a),\mathrm{T_{H}}(b)]=\mathrm{T_{H}}(\mathrm{T_{H}}(a)(b)).
\end{equation*}
It follows that $\mathcal{K}(n,\underline{t})$ is a subalgebra of
$KO(n,n+1,\underline{t}).$ Moreover, the mapping
\begin{equation*}
\rho: \mathcal{K}(n,\underline{t})\longrightarrow
HO(n,n;\underline{t}),\quad \mathrm{T_{K}}(a)\longmapsto \mathrm{T_{H}}(a)
\end{equation*}
 is an isomorphism of Lie superalgebras.
 \end{remark}
\begin{lemma}\label{l11}
Suppose $h\geq 3$. Then
\begin{eqnarray*}
KO_{h}&=&M(KO)_{h}+\sum_{ \alpha_{i}\equiv 0,1(\mathrm{mod}p)\; \forall\;i\in \overline{1,n}
\atop u\in \mathbb{B}(n),|\alpha|+|u|-2=h}\mathbb{F}\mathrm{T_{K}}(x^{(\alpha)}x^{u})\\
&&~~~~~~~~~~~+\sum_{ \alpha_{j}-\delta _{j'\in
u}\equiv0(\mathrm{mod}p)\; \forall j\in \overline{1,n},\atop u\in
\mathbb{B}(n),|\alpha|+|u|=h}\mathbb{F}\mathrm{T_{K}}(x^{(\alpha)}x^{u}x_{2n+1}).
\end{eqnarray*}
\end{lemma}
\begin{proof}
It is sufficient to show the inclusion $``\subset$. Let $\mathrm{T_{K}}(x^{(\alpha)}x^{u}x_{2n+1})$ $\in
KO_{h} \pmod{\mathcal {K}(n,\underline{t})},$ where $h\geq3,$
$\alpha\in \mathbb{A}(n,\underline{t})$ and $u\in \mathbb{B}(n).$
Note that for   $j\in \overline{1,n},$
\begin{eqnarray}\label{e3}
&&(\alpha_{j}-\delta_{j'\in
u})\mathrm{T_{K}}(x^{(\alpha)}x^{u}x_{2n+1})\nonumber\\
&=&\left[\mathrm{T_{K}}(x^{(\alpha)}x^{u}),
\mathrm{T_{K}}(x_{j}x_{j'}x_{2n+1})\right]\in M(KO)_h\pmod{\mathcal
{K}(n,\underline{t})}.
\end{eqnarray}
 If there is $j\in \overline{1,n}$
such that $\alpha_{j}-\delta_{j'\in u}\not\equiv 0\pmod{p},$ it
follows from (\ref{e3}) that
$$\mathrm{T_{K}}(x^{(\alpha)}x^{u}x_{2n+1})\in
M(KO)_h\pmod{\mathcal {K}(n,\underline{t})}.$$
Now, by Remark
\ref{r1} and Lemma \ref{l10}, one easily sees that $``\subset"$ holds.
\end{proof}

We record a general fact, which is completely analogous to the Lie algebra case (see \cite{f2,xz}). Recall that
 $L=\oplus_{i=-r}^{q}L_{i} $ denotes a finite
dimensional simple $\mathbb{Z}$-graded Lie superalgebra and $H$ a
 nilpotent subalgebra of
$L_{\bar{0}}\cap L_{0}$ with weight space decomposition
$L=\oplus_{\alpha\in \Delta}L_{(\alpha)}.$
We record a general fact, which is completely analogous to the Lie algebra case (see \cite{f2,xz}):
\begin{lemma}\label{l4}
Let $\psi: L\longrightarrow L^{\ast}$ be a $\mathbb{Z}\times
\mathrm{Map}(H,\mathbb{F})$-homogeneous skew superderivation of
degree $(l,\theta).$ Suppose
 $\psi\mid_{L^{-}}\in
\mathrm{Inn}_{\mathbb{F}}(L^{-},L^{\ast}).$ We have:

\begin{itemize}
\item [$\mathrm{(i)}$]  If $l>-q$ then $\psi$  is
inner.
\item [$\mathrm{(ii)}$] If $l=-q$ and $\Delta_{q}\cap -\Delta_{0}=\emptyset $ then
$\psi$ is inner.
 \item [$\mathrm{(iii)}$]
  If  $-2q\leq l\leq -q-1$  and
$-\Delta_{q}\not\subset\phi_{-(q+l)}$ then $\psi=0.$
\end{itemize}
\end{lemma}
We are in position to prove the main result of this paper:
\begin{theorem}\label{t2}
The second cohomology group $H^{2}(X,\mathbb{F})$ vanishes, where
$X=HO$ or $KO.$
\end{theorem}
\begin{proof}
By Theorem \ref{liu1631},  it is sufficient  to show that all the
$\mathbb{Z}$-homogeneous skew superderivations  from $X$ into
$X^{\ast}$ are inner. On the other hand, as mentioned above,
superderivationsfrom $X$ into $X^{\ast}$ of nonzero
$\mathrm{Map}(H_{X},\mathbb{F})$-degrees must be inner. Hence it is
sufficient to show  that if
 $\psi:
X\longrightarrow X^{\ast}$ is a
  skew
superderivation of  $\mathbb{Z}\times\mathrm{Map}(H_{X},\mathbb{F})$-degree $(l,\theta)$ then $\psi$ is inner.
Hence, in the below we suppose   $\psi:
X\longrightarrow X^{\ast}$ is such a superderivation.
By
Lemma \ref{l7},  $\psi\mid_{X^{-}}\in
\mathrm{Inn}_{\mathbb{F}}(X^{-},X^{\ast}).$

(1) If $l>-q,$ by Lemma \ref{l4}(i),
$\psi$ is inner.

(2) Suppose $l=-q.$ By (\ref{e5}) and (\ref{e8}), we have
$-\Delta_{X,0}\cap \Delta_{X,q}=\emptyset.$ Then by Lemma
\ref{l4}(ii),  $\psi$ is inner.

(3) Suppose $-2q\leq l\leq -q-1.$ If $l\leq-q-3,$ then $-(q+l)\geq
3.$ By Lemma \ref{l10}, for $HO$ we have
\begin{eqnarray*}
\phi_{-(q+l)}\subset\{\alpha+u&\mid& \alpha\in \mathbb{A}(n), u\in
\mathbb{B}(n), \alpha_{i}\equiv 0,1 \pmod{p},\\
&&\forall i\in \overline{1,n}, |\alpha|+|u|=-(q+l-2)\}.
\end{eqnarray*}
Combining this with (\ref{e6}), (\ref{e7}) and (\ref{e8}),
we have
$$-\Delta_{HO,q}\not\subset \phi_{-(q+l)}, \quad -\Delta_{HO,q}\not \subset \phi_{2}
\subset\Delta_{HO,2}\quad \mbox{and}\quad - \Delta_{HO,q}\not
\subset \Delta_{HO,1}=\phi_{1}. $$ By Lemma \ref{l11}, for $KO$ we
have
\begin{eqnarray*}
\phi_{-(q+l)}\subset \{\theta, \alpha+u&\mid& \alpha \in
\mathbb{A}(n), u\in \mathbb{B}(n);\;
\alpha_{i}\equiv 0,1 \pmod{p},\\
&& \forall i\in \overline{1,n}, |\alpha|+|u|=-(q+l-2)\}.
\end{eqnarray*}
Then by (\ref{e6}), (\ref{e7}) and (\ref{e8}) one easily
sees
$$-\Delta_{KO,q}\not\subset \phi_{-(q+l)}, \quad -\Delta_{KO,q}\not \subset \phi_{2}
\subset\Delta_{KO,2}\quad \mbox{and}\quad - \Delta_{KO,q}\not
\subset \Delta_{KO,1}=\phi_{1}.
 $$
By Lemma \ref{l4}(iii), $\psi=0.$ The proof is complete.
\end{proof}

\vspace{0.5cm}

\end{document}